\documentclass[11pt,a4paper]{article}

\usepackage{amsmath}
\usepackage{amsthm}
\usepackage{amssymb}
\usepackage{amscd}
\usepackage[all]{xy}

\title{On Pseudo-Effectivity of the Second Chern Classes for Smooth Threefolds
}
\author{Qihong Xie}
\date{}
\theoremstyle{plain}
\newtheorem{prop}{Proposition}[section]
\newtheorem{lem}[prop]{Lemma}
\newtheorem{thm}[prop]{Theorem}
\newtheorem{cor}[prop]{Corollary}
\newtheorem{conj}[prop]{Conjecture}
\newtheorem*{asp}{Assumption}

\theoremstyle{definition}
\newtheorem{defn}[prop]{Definition}

\newtheorem*{ack}{Acknowledgment}

\theoremstyle{remark}
\newtheorem{rem}[prop]{Remark}

\newcommand{\Q}{\mathbb Q}
\newcommand{\R}{\mathbb R}

\newcommand{\Z}{\mathbb Z}
\newcommand{\N}{\mathbb N}

\newcommand{\F}{\mathbb F}
\newcommand{\PP}{\mathbb P}
\newcommand{\OO}{\mathcal O}
\newcommand{\II}{\mathcal I}
\newcommand{\EE}{\mathcal E}
\newcommand{\NN}{\mathcal N}
\newcommand{\LL}{\mathcal L}
\newcommand{\GG}{\mathcal G}
\newcommand{\VV}{\mathcal V}
\newcommand{\TT}{\mathcal T}
\newcommand{\FF}{\mathcal F}
\newcommand{\Pic}{\mathop{\rm Pic}\nolimits}
\newcommand{\Supp}{\mathop{\rm Supp}\nolimits}

\begin{document}

\maketitle

\begin{abstract}
We prove that for smooth projective threefolds whose anticanonical divisors are nef, the second Chern classes are pseudo-effective under a weak assumption. As an application, the pseudo-effectivity of the second Chern classes implies that Kawamata's Effective Non-vanishing Conjecture holds for such threefolds.
\end{abstract}

\setcounter{section}{0}
\section{Introduction}\label{B:S1}

As is well known, Chern classes are one of the most important characteristics for complex manifolds or algebraic varieties. Therefore, it is natural to ask what kind of properties the higher Chern classes are of, if the first Chern class is assumed to satisfy some property. It is a general problem, but of great interest.

Let us consider the outcome of running the Minimal Model Program. For minimal models, as a famous result, the Bogomolov-Miyaoka-Yau inequality reveals some relations between the second and the first Chern classes (cf.\ \cite{myo}). In particular, it shows that the second Chern classes are pseudo-effective for terminal projective minimal threefolds. In general, we do not have the similar inequality and the pseudo-effectivity of the second Chern classes for Mori fiber spaces. But, if we restrict our attention to terminal projective threefolds whose anticanonical divisors are nef, then the pseudo-effectivity of the second Chern classes maybe holds.

The following Main Theorem in this paper, to some extent, gives an answer to the above problem.

\begin{thm}[Main Theorem]\label{B:1.1}
Let $X$ be a smooth projective threefold with $-K_X$ nef. Assume that ${\rm (AD_{III})}$ holds. Then the second Chern class $c_2(X)$ is pseudo-effective.
\end{thm}

Note that the Main Theorem holds only when the following assumption ${\rm (AD_{III})}$ is true. As a special case, we prove that ${\rm (AD_{III})}$ holds if $\rho(X)\leq 3$.

\begin{asp}[${\bf AD_{III}}$]\label{B:1.2}
Let $X$ be a smooth projective threefold with $-K_X$ nef, $f:X\rightarrow Y$ an extremal contraction with the extremal ray $R=\R_+[l]$, which contracts a divisor $E$ to a curve $C$ such that either

{\rm (A)} $C\cong \PP^1, \quad \NN_{C|Y}\cong \OO(-1)\oplus \OO(-2),$ or

{\rm (B)} $C\cong \PP^1, \quad \NN_{C|Y}\cong \OO(-2)\oplus \OO(-2).$

Then there exists a positive integer $n$ such that $c_2(X)+nl$ is pseudo-effective.
\end{asp}

We say that an extremal contraction $f: X\rightarrow Y$ is good, if we can prove that $c_2(X)$ is pseudo-effective, or at least, that there exists a positive integer $n$ such that $c_2(X)+nl$ is pseudo-effective, where $\R_+[l]$ is the corresponding extremal ray of $f$. The main idea of the proof of Theorem \ref{B:1.1} is to investigate the goodness of all extremal contractions from $X$, which implies the pseudo-effectivity of $c_2(X)$.

In \S \ref{B:S2}, we will give some necessary definitions and propositions, and obtain an application to the Effective Non-vanishing Conjecture. In \S \ref{B:S3} and \S \ref{B:S4}, we will prove Theorem \ref{B:1.1} when the irregularity $q(X)=1,0$, respectively.

In the whole paper, we will use freely the results on the Minimal Model Theory. We refer to \cite{kmm} and \cite{km} for the details.

We work over the field of complex numbers.

\begin{ack}
I would like to express my gratitude to Professor Yujiro Kawamata for his valuable advice and warm encouragement. I also thank Professor Keiji Oguiso and Dr.\ Yasunari Nagai for stimulating discussions. In addition, I am very grateful to the referee for his helpful suggestions and generous linguistic advice.
\end{ack}

\section{Preliminaries and application}\label{B:S2}

\begin{defn}\label{B:2.1}
Let $X$ be a proper variety. A 1-cycle is a formal linear combination of irreducible, reduced and proper curves $C=\sum a_iC_i$. A 1-cycle is said to be effective if $a_i\geq 0$ for every $i$. Two 1-cycles $C,C'$ are said to be numerically equivalent if $C.D=C'.D$ for any Cartier divisor $D$. We denote by $N_1(X)$ the $\R$-vector space generated by all 1-cycles with real coefficients modulo numerical equivalence. Let $NE(X)$ be the convex subcone of $N_1(X)$ generated by all effective 1-cycles, and $\overline{NE}(X)$ the closure of $NE(X)$ in $N_1(X)$.

A 1-cycle $C$ is said to be pseudo-effective, if its numerical equivalence class is contained in $\overline{NE}(X)$. If $X$ is projective, then a 1-cycle $C$ is pseudo-effective if and only if $C.H\geq 0$ for any ample Cartier divisor $H$ on $X$ by Kleiman's Ampleness Criterion (cf.\ \cite{kl}).

A $\Q$-Cartier divisor $D$ is said to be nef, if $D.C\geq 0$ for any irreducible curve $C$. For any nef $\Q$-Cartier divisor $D$ on $X$, the numerical dimension $\nu(D)$ is defined to be the greatest non-negative integer $\nu$ such that $D^\nu\not\equiv 0$.

$X$ is called a terminal variety, if $X$ has only terminal singularities.
\end{defn}

\begin{thm}\label{B:2.2}
Let $X$ be a terminal projective threefold with $-K_X$ nef. Then the following conclusions hold.

{\rm (i)} If $\nu(-K_X)=0$, then $c_2(X)$ is pseudo-effective {\rm (cf.\ \cite{myo})};

{\rm (ii)} If $\nu(-K_X)=1$, then $c_2(X)$ is pseudo-effective {\rm (cf.\ \cite{matsuki})};

{\rm (iii)} If $\nu(-K_X)=3$, then $c_2(X)$ is pseudo-effective {\rm (cf.\ \cite{kmmt})};

{\rm (iv)} $c_1(X).c_2(X)\geq 0$, hence $\chi(\OO_X)\geq 0$ {\rm (cf.\ \cite{kmmt}, \cite{matsuki})}.
\end{thm}

\begin{proof}
In the above papers, the pseudo-effectivity of $c_2(X)$ has not been explicitly mentioned, but it is easy to derive these conclusions.
\end{proof}

Theorem \ref{B:2.2} enables us to put forward the following conjecture.

\begin{conj}\label{B:2.3}
Let $X$ be a terminal projective threefold with $-K_X$ nef. Then the second Chern class $c_2(X)$ is pseudo-effective.
\end{conj}

As a partial answer, we will prove in Theorem \ref{B:1.1} that Conjecture \ref{B:2.3} holds in the smooth case under a weak assumption.

It follows from Theorem \ref{B:2.2} that we only need to verify the case when $\nu(-K_X)=2$ for proving Conjecture \ref{B:2.3}. First, we divide this case into more explicit subcases by the irregularity $q(X):=\dim H^1(X,\OO_X)$.

\begin{prop}\label{B:2.4}
Let $X$ be a terminal projective threefold such that $-K_X$ is nef and $\nu(-K_X)=2$. Then either

{\rm (i)} $q(X)=1$ and $X$ is Gorenstein, or

{\rm (ii)} $q(X)=0$.
\end{prop}

\begin{proof}
Let $H$ be an ample Cartier divisor on $X$. We consider the following exact sequence:
\[ 0\rightarrow\OO_X(mK_X-H)\rightarrow\OO_X(mK_X)\rightarrow\OO_H(mK_X\vert_H)\rightarrow 0\]
where $m$ is any integer.

Since $-(mK_X-H)=-mK_X+H$ is ample for $m\geq 0$, $-mK_X\vert_H$ is nef and big for $m\geq 1$, it follows from the Kawamata-Viehweg vanishing theorem that
$H^i(\OO_X(mK_X-H))=0$ for $i<3$ and $m\geq 0$, and
$H^i(\OO_H(mK_X\vert_H))=0$ for $i<2$ and $m\geq 1$.
Therefore it follows from the above exact sequence that $H^0(mK_X)=H^1(mK_X)=0$ for $m\geq 1$, namely, $h^i(-mK_X)=0$ for $i=2,3$ and $m\geq 0$.

By Theorem \ref{B:2.2}(iv), we have that $\chi(\OO_X)=\sum (-1)^ih^i(\OO_X)=h^0(\OO_X)-h^1(\OO_X)=1-q(X)\geq 0$. If $q(X)=1$, then $\chi(\OO_X)=0$, which implies that $X$ is Gorenstein (cf.\ \cite{matsuki}), and there is the Albanese map $\alpha: X\rightarrow {\rm Alb}(X)$ to an elliptic curve. Otherwise $q(X)=0$, which completes the proof.
\end{proof}

The pseudo-effectivity of the second Chern classes is closely related to the Effective Non-vanishing Conjecture, which has been put forward formally by Yujiro Kawamata (cf.\ \cite{ka00}).

\begin{conj}\label{B:2.5}
Let $X$ be a complete normal variety, $B$ an effective $\R$-divisor on $X$ such  that the pair $(X,B)$ is Kawamata log terminal (KLT, for short), and $D$ a Cartier divisor on $X$. Assume that $D$ is nef and that $D-(K_X+B)$ is nef and big. Then $H^0(X,D)\neq 0$.
\end{conj}

Kawamata has proven that the Effective Non-vanishing Conjecture holds for all log surfaces with only KLT singularities. As an application, we prove that the pseudo-effectivity of the second Chern classes implies that the Effective Non-vanishing Conjecture holds for terminal projective threefolds whose anticanonical divisors are nef.

\begin{prop}\label{B:2.6}
Let $X$ be a terminal projective threefold with $-K_X$ nef, $D$ a Cartier divisor on $X$. Assume that $D$ is nef and $D-K_X$ is big. Then the pseudo-effectivity of $c_2(X)$ implies that $H^0(X,\OO_X(D))\neq 0$.
\end{prop}

\begin{proof}
By the Kawamata-Viehweg vanishing theorem, we have $H^i(X,D)=0$ for any positive integer $i$. Thus the condition $H^0(X,D)\neq 0$ is equivalent to saying that $\chi(X,D)\neq 0$.

If $\nu(D)<3$, then we may reduce this case to the log surface case which has been proven by Kawamata. Assume that $D$ is nef and big. It follows from the Riemann-Roch Theorem and the pseudo-effectivity of $c_2(X)$ that
\[ h^0(\OO_X(D))=\frac {1} {12}D(D-K_X)(2D-K_X)+\frac {1} {12}D.c_2(X)+\chi(\OO_X)>0. \]
\end{proof}

\section{Proof of the case $q(X)=1$}\label{B:S3}

\begin{defn}\label{B:3.1}
Let $X$ be a variety, $\pi: X\rightarrow A$ a surjective morphism to a curve $A$. A curve $C\subset X$ is said to be an \'etale multi-section of $\pi$, if $\pi|_C:C\rightarrow A$ is a finite \'etale cover.
\end{defn}

In fact, in the case $q(X)=1$, the structure of $X$ is determined by the following theorem (cf.\ \cite{ps}, Corollary 3.4).

\begin{thm}\label{B:3.2}
Let $X$ be a smooth projective threefold such that $-K_X$ is nef, $\nu(-K_X)=2$ and $q(X)=1$. Let $\alpha: X\rightarrow {\rm Alb}(X)=A$ be the Albanese map to a smooth elliptic curve $A$. Then there exists a sequence of blow-ups $\varphi_i:X_i\rightarrow X_{i+1}$, $0\leq i\leq s$, with $X_0=X$ and inducing morphisms $\alpha_i: X_i\rightarrow A$, such that

{\rm (D)} all $X_i$ are smooth with $-K_{X_i}$ nef, and $\varphi_i$ is the blow-up of a smooth curve $C_i$, which is an \'etale multi-section of $\alpha_{i+1}: X_{i+1}\rightarrow A$.

{\rm (F)} the induced morphism $\alpha_{s+1}: X_{s+1}\rightarrow A$ is one of the following cases:

{\rm (I)} a $\PP^2$-bundle;

{\rm (II)} a $\PP^1\times\PP^1$-bundle;

{\rm (III)} $\alpha_{s+1}$ factors as $h\circ g$ with $g: X_{s+1}\rightarrow Y$ a conic bundle and $h: Y\rightarrow A$ a $\PP^1$-bundle.
\end{thm}

\begin{proof}
It follows from the proofs of Proposition 1.7 and Theorem 3.3 of \cite{ps} that we have to exclude the following case: 

(H) $f: X\rightarrow S$ is a $\PP^1$-bundle, and $S$ is a hyperelliptic surface.

We prove that case (H) cannot occur. Otherwise, we may take a finite \'etale cover $\pi: W\rightarrow S$ such that $W$ is an abelian surface. Let $Y=X\times_S W$ be the fiber product over $S$, $g: Y\rightarrow W$ the induced morphism. It is easy to show that $-K_Y$ is nef and $\nu(-K_Y)=2$. It follows from Proposition \ref{B:2.4} and Hodge symmetry that $h^0(\Omega^2_Y)=h^2(\OO_Y)=0$. Since $h^0(K_W)=1$, there is a nowhere vanishing 2-form on $W$. Hence $Y$ has a nonzero 2-form by pullback with $g$, this is absurd.
\end{proof}

The keypoint of the proof of the pseudo-effectivity of $c_2(X)$ is a direct verification for $X_{s+1}$ and using induction on $i$ for the general case.

\begin{lem}[Case F-I]\label{B:3.3}
Let $C$ be a smooth elliptic curve, $\EE$ a locally free sheaf of rank 3 on $C$. Assume that the $\PP^2$-bundle $X=\PP_C(\EE)$ is a smooth threefold with $-K_X$ nef. Then $c_2(X)$ is pseudo-effective.
\end{lem}

\begin{proof}
Let $L$ be the divisor corresponding to the tautological line bundle $\OO_X(1)$, $F$ the fiber of $\pi: X\rightarrow C$, and $r=\deg\EE=\deg c_1(\EE)$.

We have the following exact sequences:
\begin{eqnarray}
0\rightarrow\Omega_{X/C}\rightarrow\pi^*\EE\otimes\OO_X(-1)\rightarrow\OO_X\rightarrow 0 \label{B:eqn1} \\
0\rightarrow \pi^*\Omega_C\rightarrow\Omega_X\rightarrow\Omega_{X/C}\rightarrow 0 \label{B:eqn2}
\end{eqnarray}
It follows from (\ref{B:eqn2}) that
\begin{eqnarray}
c_1(\Omega_X) & = & \pi^*c_1(\Omega_C)+c_1(\Omega_{X/C})=c_1(\Omega_{X/C}) \nonumber \\
c_2(\Omega_X) & = & \pi^*c_1(\Omega_C).c_1(\Omega_{X/C})+c_2(\Omega_{X/C})=c_2(\Omega_{X/C}) \nonumber
\end{eqnarray}
It follows from (\ref{B:eqn1}) that
\begin{eqnarray}
c_1(\Omega_{X/C}) & = & c_1(\pi^*\EE\otimes\OO_X(-1))=rF-3L \nonumber \\
c_2(\Omega_{X/C}) & = & c_2(\pi^*\EE\otimes\OO_X(-1))=-2rF.L+3L^2=\frac{1}{3}c_1^2(\Omega_{X/C}) \nonumber
\end{eqnarray}
Hence $c_2(X)=(-K_X)^2/3$ is pseudo-effective since $-K_X$ is nef.
\end{proof}

\begin{rem}\label{B:3.4}
At first, such a subcase does exist. For example, let $\EE=\OO_C^{\oplus 3}$, $X\cong C\times\PP^2$. It is easy to show that $-K_X$ is nef and $\nu(-K_X)=2$. Secondly, given a multi-polarization ($H_1,H_2$) with $H_1,H_2$ ample divisors on $X$. Then in Lemma 3.3, the tangent bundle $\TT_X$ is always unstable with a destabilizing subsheaf $\TT_{X/C}$.
\end{rem}

\begin{lem}[Case F-II]\label{B:3.5}
Let $C$ be a smooth elliptic curve, $\pi: X\rightarrow C$ a $\PP^1\times\PP^1$-bundle. Assume that $-K_X$ is nef. Then $c_2(X)$ is pseudo-effective.
\end{lem}

\begin{proof}
By the definition of $\pi$, there exists a locally free sheaf $\EE$ of rank 4, such that $Y=\PP_C(\EE)$ is a $\PP^3$-bundle. $X\subset Y$ is a divisor on $Y$ such that $X_p\subset Y_p\cong\PP^3$ is isomorphic to $\PP^1\times\PP^1$ over each point $p\in C$. Let $L$ be the divisor corresponding to the tautological line bundle $\OO_Y(1)$, $F$ the fiber of $\pi: Y\rightarrow C$, and $r=\deg\EE=\deg c_1(\EE)$. Then $X\sim 2L$.

We have similar exact sequences to (\ref{B:eqn1}) and (\ref{B:eqn2}) for $Y$. Then
\begin{eqnarray}
c_1(\Omega_Y) & = & c_1(\Omega_{Y/C})=rF-4L \nonumber \\
c_2(\Omega_Y) & = & c_2(\Omega_{Y/C})=-3rF.L+6L^2 \nonumber
\end{eqnarray}
We also have the following exact sequence:
\begin{eqnarray}
0\rightarrow \NN^*_{X|Y}\rightarrow \Omega_Y\otimes\OO_X \rightarrow \Omega_X \rightarrow 0 \label{B:eqn3}
\end{eqnarray}
There are some simple computations from (\ref{B:eqn3}):
\begin{eqnarray}
c_1(\Omega_X) & = & K_X=(K_Y+X)|_X=(rF-2L)|_X \nonumber \nonumber \\
c_2(\Omega_X) & = & c_2(\Omega_Y\otimes\OO_X)-c_1(\NN^*_{X|Y}).c_1(\Omega_X) \nonumber \\
              & = & c_2(\Omega_Y)|_X+X|_X.K_X \nonumber \\
              & = & (-rF.L+2L^2)|_X \nonumber \\
              & = & (-K_X).L|_X \nonumber
\end{eqnarray}
Hence $c_2(X)$ is pseudo-effective since $-K_X$ is nef.
\end{proof}

\begin{lem}[Case F-III]\label{B:3.6}
Let $C$ be a smooth elliptic curve, $S$ a $\PP^1$-bundle over $C$. Let $f:X\rightarrow S$ be a conic bundle, $X$ a smooth threefold with $-K_X$ nef. Then $c_2(X)$ is pseudo-effective.
\end{lem}

\begin{proof}
Let $\Delta\subset S$ be the discriminant locus of $f$. Then we have the following exact sequence (cf.\ \cite{st}):
\begin{eqnarray}
0 \rightarrow f^*\Omega_S \rightarrow \Omega_X \rightarrow \OO_X(K_{X/S}) \rightarrow \OO_\Gamma \rightarrow 0 \label{B:eqn4}
\end{eqnarray}
where $K_{X/S}=K_X-f^*K_S$ is the relative canonical divisor, $\Gamma$ is a locally complete closed subscheme of $X$ of pure dimension 1 with $f(\Gamma)=\Delta$. The restriction $f|_{\Gamma\setminus f^{-1}(\Delta_{\rm sing})}:\Gamma\setminus f^{-1}(\Delta_{\rm sing}) \rightarrow \Delta_{\rm reg}$ is an isomorphism and $\Gamma\cap X_s=(X_s)_{\rm red}$ for all $s\in\Delta_{\rm sing}$.

It is easy to see that $c_2(\Omega_S)=K_S^2=0$. Furthermore, we have that $-(4K_S+\Delta)$ is nef by \cite{dps}. It follows from (\ref{B:eqn4}) and Lemma \ref{B:3.7} that
\begin{eqnarray}
c_2(\Omega_X) & = & f^*c_2(\Omega_S)+f^*c_1(\Omega_S).K_{X/S}-c_2(\OO_\Gamma) \nonumber \\
              & = & f^*(-K_S).(-K_X)+\Gamma \nonumber \\
              & = & f^*(-K_S-\frac{1}{4}\Delta).(-K_X)+\frac{1}{4}f^*\Delta.(-K_X)+\Gamma. \nonumber
\end{eqnarray}
Thus $c_2(X)$ is pseudo-effective.
\end{proof}

\begin{lem}\label{B:3.7}
Let $X$ be a smooth projective threefold, $\Gamma$ a locally complete closed subscheme of X of pure dimension 1. Then $c_2(\OO_\Gamma)=-\Gamma$.
\end{lem}

\begin{proof}
Both sides are clearly additive over subschemes with disjoint supports. If $\gamma$ is a smooth curve in $X$, then there exist smooth hypersurface sections $H_1,H_2$ such that $\gamma\subset\Supp(H_1\cap H_2)$, every irreducible component of $H_1\cap H_2$ is smooth in $X$ and $H_1,H_2$ meet transversally. Let $Y=H_1\cap H_2$. From the exact sequence
\begin{eqnarray}
0 \rightarrow \OO(-H_1-H_2) \rightarrow \OO(-H_1)\oplus\OO(-H_2) \rightarrow \II_Y \rightarrow 0 \label{B:eqn5}
\end{eqnarray}
we can calculate that $c_2(\II_Y)=H_1.H_2$, hence $c_2(\OO_Y)=-H_1.H_2$. Since $c_2$ is invariant in an algebraic family, all irreducible curves in $Y$ are algebraically equivalent on $H_1$, therefore $c_2(\OO_\gamma)=-\gamma$.

In the general case, since $\Gamma$ is a local complete intersection, there exists a sufficiently ample divisor $H$, such that $\OO_\Gamma(H)$ is generated by global sections and there exist $H_1,H_2\in |H|$ whose local equations generate the ideal of $\Gamma$ in $\OO_{X,\gamma}$ for each irreducible curve $\gamma\in\Gamma$, and all other intersections are transversal. Let $Y$ be the scheme theoretic intersection of $H_1$ and $H_2$. Then from the exact sequence (\ref{B:eqn5}), we have $c_2(\OO_Y)=-H_1.H_2$ as before, and each irreducible curve $l$ of $\Supp(Y)\setminus\Supp(\Gamma)$ contributes $-l$. Thus $c_2(\OO_\Gamma)=-\Gamma$ in the general case.
\end{proof}

\begin{prop}\label{B:3.8}
Let $Y$ be a smooth projective threefold. Let $f:X\rightarrow Y$ be the blow-up along a nonsingular subvariety $Z\subset Y$, $E$ the exceptional divisor of $f$ on $X$. Then the following assertions hold.

{\rm (i)} If $Z=p$ is a point, then $c_2(X)=f^*c_2(Y)$;

{\rm (ii)} If $Z=C$ is a curve, then $c_2(X)=f^*c_2(Y)-E^2-\deg c_1(C)F$, $f^*C=-E^2+\deg c_1(\NN_{C|Y})F$, where $F$ is the fiber of $f|_E:E\rightarrow C$.
\end{prop}

\begin{proof}
This follows from \cite{fu} and \cite{myn}.
\end{proof}

The proof of Proposition \ref{B:3.9} is almost identical to that of Proposition 3.3 of \cite{dps}.

\begin{prop}\label{B:3.9}
Let $Y$ be a smooth projective threefold, $f:X\rightarrow Y$ be the blow-up along a smooth curve $C$ in $Y$ such that $-K_X$ is nef. Let $H$ be a nef divisor on $X$. Then $f_*H$ is nef possibly except the following cases:

{\rm (A)} $C\cong \PP^1, \quad \NN_{C|Y}\cong \OO(-1)\oplus \OO(-2);$

{\rm (B)} $C\cong \PP^1, \quad \NN_{C|Y}\cong \OO(-2)\oplus \OO(-2);$

{\rm (C)} $C\cong \PP^1, \quad \NN_{C|Y}\cong \OO(-1)\oplus \OO(-1).$
\end{prop}

\begin{proof}
Let $E$ be the exceptional divisor of $f$. Then $K_X=f^*K_Y+E$. Let $\NN=\NN_{C|Y}$ be the normal bundle of $C$ in $Y$, $\VV=\NN^*\otimes \LL$ with $\LL\in\Pic(C)$ its normalization, i.e. $H^0(\VV)\neq 0, H^0(\VV\otimes\GG)=0$ for all $\GG\in\Pic(C)$ with $\deg\GG<0$. Then $E=\PP(\NN^*)\cong\PP(\VV)$, and the tautological line bundle $\OO_{\PP(\VV)}(1)$ has a section $C_1$ such that $C_1^2=-e=c_1(\VV)$. Let $\mu=\deg\LL$, $F$ the fiber of the $\PP^1$-bundle $f|_E:E\rightarrow C$.

Let $-K_X|_E\equiv aC_1+bF$, since $K_X.F=-1$, we have $a=1$. Moreover, $\NN_{E|X}=E|_E\equiv -C_1+\mu F$ from $\OO_E(E|_E)\cong \OO_E(-1)$ and the definition of $\mu$. Let $g$ be the genus of $C$, we have $K_E^2=8(1-g)$. On the other hand, by the adjunction formula, $K_E^2=((K_X+E)|_E)^2=(-2C_1+(\mu-b)F)^2$. Hence
\begin{eqnarray}
e+\mu-b & = & 2(g-1) \label{B:eqn6}
\end{eqnarray}
Since $-K_X|_E$ is nef on $E$, $-K_X|_E.C_1\geq 0$ implies that
\begin{eqnarray}
b & \geq & e \label{B:eqn7}
\end{eqnarray}
and $(-K_X|_E)^2\geq 0$ implies that
\begin{eqnarray}
b & \geq & \frac{1}{2}e \label{B:eqn8}
\end{eqnarray}

Let $C'$ be an irreducible curve in $Y$. If $C'\neq C$, then $f_*H.C'=H.f^*C'\geq H.f_*^{-1}C'\geq 0$ by the projection formula, where $f_*^{-1}C'$ is the strict transform of $C'$. Otherwise,
\begin{eqnarray}
f_*H.C & = & H.f^*C=H.(-E^2+\deg c_1(\NN_{C|Y})F) \nonumber \\
       & = & H.(C_1-\mu F+(e+2\mu)F) \nonumber \\
       & = & H.(C_1+(b+2(g-1))F)  \label{B:eqn9}
\end{eqnarray}

Note that if $\VV$ is decomposable then $e\geq 0$, otherwise $e\geq -2g$. Anyway, we have that $b\geq e/2\geq -g$, hence $b+2(g-1)\geq g-2$.

(\ref{B:3.9}.1) If $g\geq 2$, then $f_*H.C\geq 0$;

(\ref{B:3.9}.2) If $g=1$, then $e\geq -1$, hence $b\geq e/2\geq -1/2$, i.e. $b\geq 0$. $f_*H.C=H.(C_1+bF)\geq 0$;

(\ref{B:3.9}.3) If $g=0$ and $b\geq 2$, then $f_*H.C=H.(C_1+(b-2)F)\geq 0$. Since $C$ is rational, we may assume that $0\leq e\leq b\leq 1$, and obtain the cases (A), (B) and (C).
\end{proof}

\begin{prop}\label{B:3.10}
With the notation as above. Assume that both $-K_X$ and $-K_Y$ are nef but not big. Then $c_2(X)=f^*c_2(Y)+C_1$.
\end{prop}

\begin{proof}
It follows from Proposition \ref{B:3.8} that
\begin{eqnarray}
c_2(X) & = & f^*c_2(Y)-E^2-\deg c_1(C)F \nonumber \\
       & = & f^*c_2(Y)+(C_1-\mu F)+(2g-2)F \nonumber \\
       & = & f^*c_2(Y)+C_1+(e-b)F \nonumber
\end{eqnarray}
The formula $(-K_X)^3=(-K_Y)^3-2((-K_Y).C-g(C)+1)$ implies that $(-K_Y).C=g-1\geq 0$, since both $-K_X$ and $-K_Y$ are nef but not big. On the other hand, $(-K_Y).C=(-f^*K_Y).C_1=(-K_X+E).C_1=b+\mu$. It follows from (\ref{B:eqn6}) and (\ref{B:eqn8}) that $2b-e=1-g\geq 0$. Therefore $g=1$, $2b=e\geq -1$, hence $b\geq 0,e\geq 0$. It follows from (\ref{B:eqn7}) that $b=e=0$, which completes the proof.
\end{proof}

\begin{rem}\label{B:3.11}
In Propositions \ref{B:3.9} and \ref{B:3.10}, we need not assume that $q(X)>0$. If $q(X)=1$, then the conclusion $g=1$ in Proposition \ref{B:3.10} also follows from Proposition 3.2 of \cite{ps}.
\end{rem}

\begin{thm}\label{B:3.12}
Let $X$ be a smooth projective threefold such that $-K_X$ is nef, $\nu(-K_X)=2$ and $q(X)=1$. Then $c_2(X)$ is pseudo-effective.
\end{thm}

\begin{proof}
It is necessary to verify case (D) in Theorem \ref{B:3.2}. Let $\varphi_i: X_i\rightarrow X_{i+1}$ be a blow-up along an elliptic curve $C_i$. Note that $(-K_{X_i})^3=(-K_{X_{i+1}})^3=0$ since otherwise $q(X_i)=0$ resp.\ $q(X_{i+1})=0$ by the Kawamata-Viehweg vanishing theorem. It follows from Proposition \ref{B:3.10} that $c_2(X_i)=\varphi_i^* c_2(X_{i+1})+C^1_i$, where $C^1_i$ is the canonical section of $C_i$ for $\varphi_i$. Let $H$ be any nef divisor on $X_i$, then $(\varphi_i)_*H$ is nef by Proposition \ref{B:3.9}. Assume that $c_2(X_{i+1})$ is pseudo-effective, then $c_2(X_i).H=c_2(X_{i+1}).(\varphi_i)_*H+C^1_i.H\geq 0$, namely $c_2(X_i)$ is pseudo-effective. Induction on $i$ completes the proof.
\end{proof}

\section{Proof of the case $q(X)=0$}\label{B:S4}

The case $q(X)=0$ is more complicated than the case $q(X)=1$, because, at least, we cannot give a nice classification for such $X$. But we may take an extremal contraction from $X$, and investigate its goodness. Next, we will make use of the general theory of extremal contractions from smooth projective threefolds given by Shigefumi Mori (cf.\ \cite{mo}).

\begin{prop}\label{B:4.1}
Let $X$ be a smooth projective threefold such that $-K_X$ is nef, $\nu(-K_X)=2$ and $q(X)=0$. Let $f:X\rightarrow Y$ be an extremal contraction. Then $(Y,f)$ is one of the following cases:

${\rm (F_I)}$ $f$ is a del Pezzo fibration, and $Y\cong\PP^1$;

${\rm (F_{II})}$ $f$ is a conic bundle with discriminant locus $\Delta$ (possibly empty), such that $-(4K_Y+\Delta)$ is nef, and $Y\cong\PP^2$;

${\rm (F_{III})}$ $f$ is a conic bundle with discriminant locus $\Delta$ (possibly empty), such that $-(4K_Y+\Delta)$ is nef, and there is a morphism $\alpha:Y\rightarrow\F_n$, which is the composition of a sequence of one point blow-ups over the Hirzebruch surface $\F_n$;

${\rm (D_I)}$ $f$ is a birational morphism which contracts an irreducible divisor $E$ to a point $p\in Y$, and $Y$ is a terminal projective threefold with $-K_Y$ nef and big;

${\rm (D_{II})}$ $f$ is a birational morphism which contracts an irreducible divisor $E$ to a smooth curve $C\subset Y$, and $Y$ is a smooth projective threefold with $-K_Y$ nef;

${\rm (D_{III})}$ $f$ is a birational morphism which contracts an irreducible divisor $E$ to a smooth curve $C\subset Y$, and $Y$ is a smooth projective threefold, such that $-K_Y$ is nef except along $C$. Furthermore, $C$ is one of the following two cases:

{\rm (A)} $C\cong \PP^1, \quad \NN_{C|Y}\cong \OO(-1)\oplus \OO(-2);$

{\rm (B)} $C\cong \PP^1, \quad \NN_{C|Y}\cong \OO(-2)\oplus \OO(-2).$
\end{prop}

\begin{proof}
It follows from the definition and Proposition \ref{B:2.4} that $h^i(\OO_X)=0$ for any integer $i>0$. Since $f$ is extremal, $R^jf_*\OO_X=0$ for any integer $j>0$ by the Kawamata-Viehweg vanishing theorem. The standard argument by making use of the Leray spectral sequence deduces that $h^i(\OO_Y)=0$ for any integer $i>0$. 

(\ref{B:4.1}.0) It is easy to see that $\dim Y>0$, since $X$ cannot be a Fano threefold.

(\ref{B:4.1}.1) If $\dim Y=1$, then $f$ is a del Pezzo fibration, and $Y$ is a smooth curve. $h^1(\OO_Y)=0$ implies that $Y\cong\PP^1$.

(\ref{B:4.1}.2) If $\dim Y=2$, then $f$ is a conic bundle with discriminant locus $\Delta$ (possibly empty), and $Y$ is a smooth surface.

First, we prove that $\kappa(Y)=-\infty$.

If $\Delta\neq\emptyset$, then we have $-(4K_Y+\Delta)$ is nef (cf.\ \cite{dps}). We claim that $\kappa(Y)=-\infty$. Otherwise, there exists an ample divisor $H$ on $Y$, such that $K_Y.H\geq 0$. But $-(4K_Y+\Delta).H\geq 0$ implies that $K_Y.H\leq -(1/4)\Delta.H<0$. This is absurd.

If $\Delta=\emptyset$, then we have $-K_Y$ is nef (cf.\ \cite{dps}). We may contract all $(-1)$-curves to get a birational morphism $g:Y\rightarrow Z$. If $\kappa(Y)\geq 0$, then $K_Z$ is nef by definition. On the other hand, it is easy to see that $-K_Z$ is nef. Hence $K_Z\equiv 0$, and $Y\cong Z$ since $K_Y^2=K_Z^2=0$. It follows from $h^1(\OO_Y)=h^2(\OO_Y)=0$ that $q(Y)=h^0(K_Y)=0$, hence $Y$ is an Enriques surface by the Classification Theory of surfaces. We prove that this case cannot occur. Let $\pi: S\rightarrow Y$ be a degree 2 \'etale cover from a $K$3 surface $S$ to $Y$, $h: W=X\times_Y S\rightarrow S$ the fiber product over $Y$. Note that the projection $\tau:W\rightarrow X$ is \'etale since $\pi$ is \'etale. Thus $-K_W$ is nef and $\nu(-K_W)=2$. It follows from Proposition \ref{B:2.4} and Hodge symmetry that $h^0(W,\Omega_W^2)=h^2(\OO_W)=0$. On the other hand, $S$ has a nowhere vanishing 2-form since $S$ is $K$3, then $W$ has a nonzero 2-form by pullback with $h$. This is absurd.

Therefore $p_2(Y)=q(Y)=0$ which implies that $Y$ is a rational surface by Castelnuovo's Rationality Criterion. Thus there is a morphism $\alpha:Y\rightarrow \F_n$, which is the composition of a sequence of one point blow-ups over the Hirzebruch surface $\F_n$, or $Y\cong\PP^2$.

(\ref{B:4.1}.3) If $\dim Y=3$, then $f$ is a divisorial contraction.

(\ref{B:4.1}.3.1) If $\dim f(E)=0$, then we can write that $K_X=f^*K_Y+aE$, where $a$ is a positive rational number. It is easy to verify that $-K_Y$ is nef. $K_X^3=K_Y^3+a^3E^3$ implies that $(-K_Y)^3=a^3(-E|_E)^2>0$, hence $-K_Y$ is nef and big. The value of $a$ and the explicit structure of $p=f(E)\in Y$ can be found in \cite{mo}.

(\ref{B:4.1}.3.2) If $\dim f(E)=1$, then $f$ is just the blow-up of $Y$ along a smooth curve $C=f(E)$. It follows from \cite{dps} that $-K_Y$ is nef except for the two cases (A) and (B) listed in ${\rm (D_{III})}$.
\end{proof}

From now on, we assume that $\R_+[l]$ is the extremal ray which induces the extremal contraction $f: X\rightarrow Y$ given in Proposition \ref{B:4.1}.

\begin{lem}\label{B:4.2}
In case ${\rm (F_I)}$, let $f:X\rightarrow\PP^1$ be a del Pezzo fibration. Then $c_2(X)$ is pseudo-effective.
\end{lem}

\begin{proof}
Since $\Pic(X)=f^*\Pic(\PP^1)\oplus\Z$ and $-K_X$ is $f$-ample, then for any ample divisor $M$ on $X$, we have 
\[ M \equiv a(-K_X)+bX_\xi \]
where $a,b\in\Q$, and $X_\xi=f^{-1}(\xi)$ for a general point $\xi\in\PP^1$. Then $M.l=a(-K_X).l>0$ implies that $a>0$, and $M.(-K_X)^2=bX_\xi.(-K_X)^2=b(-K_{X_\xi})^2>0$ implies that $b>0$, since $(-K_X.l)>0$ and $X_\xi$ is a smooth del Pezzo surface.

It is sufficient to prove that $c_1(X).c_2(X)\geq 0$ and $X_\xi.c_2(X)\geq 0$ for proving the pseudo-effectivity of $c_2(X)$. It is obvious that $c_1(X).c_2(X)=24\chi(\OO_X)=24>0$. Let $i:X_\xi\rightarrow X$ be the closed immersion. Since $X_\xi$ is smooth, there is an exact sequence:
\begin{eqnarray}
0 \rightarrow \TT_{X_\xi} \rightarrow \TT_X\otimes\OO_{X_\xi} \rightarrow \NN_{X_\xi|X} \rightarrow 0 \label{B:eqn10}
\end{eqnarray}
It follows from (\ref{B:eqn10}) that
\begin{eqnarray}
X_\xi.c_2(X) & = & i^*c_2(X)=c_2(\TT_X\otimes\OO_{X_\xi}) \nonumber \\
             & = & c_2(\TT_{X_\xi})+c_1(\TT_{X_\xi}).c_1(\NN_{X_\xi|X}) \nonumber \\
             & = & c_2(\TT_{X_\xi})\geq 3. \nonumber
\end{eqnarray}
\end{proof}

\begin{lem}\label{B:4.3}
In case ${\rm (F_{II})}$, let $f:X\rightarrow \PP^2$ be a conic bundle. Then $c_2(X)$ is pseudo-effective.
\end{lem}

\begin{proof}
Since $\Pic(X)=f^*\Pic(\PP^2)\oplus\Z$ and $-K_X$ is $f$-ample, then for any ample divisor $M$ on $X$, we have 
\[ M \equiv a(-K_X)+bF \]
where $a,b\in\Q$, and $F=f^{-1}(H)$, where $H$ is a line in $\PP^2$. For the extremal ray $\R_+[l]$, we assume that $l$ is a smooth conic. Then $(M.l)=a(-K_X.l)>0$ implies that $a>0$, and $M.(-K_X)^2=bF.(-K_X)^2=b(-K_X|_F)^2>0$ implies that $b>0$, since $(-K_X.l)>0$ and $-K_X|_F$ is nef on $F$.

It is sufficient to prove that $c_1(X).c_2(X)\geq 0$ and $F.c_2(X)\geq 0$ for proving the pseudo-effectivity of $c_2(X)$. It is obvious that $c_1(X).c_2(X)=24\chi(\OO_X)=24>0$. It follows from Lemma \ref{B:3.6} that the following equality holds for any conic bundle $f:X\rightarrow Y$.
\begin{eqnarray}
c_2(X)=f^*(c_2(Y)-c_1^2(Y))+f^*(-K_Y).(-K_X)+\Gamma \label{B:eqn11}
\end{eqnarray}
In this case, $Y=\PP^2$, hence $c_2(X)=-6l+3F.(-K_X)+\Gamma$. Since $H$ is very ample, we may assume that $H$ intersects $\Delta$ transversally. Then $F.\Gamma=H.\Delta$, and $F.c_2(X)=3l.(-K_X)+F.\Gamma=6+H.\Delta>0$.
\end{proof}

\begin{lem}\label{B:4.4}
In case ${\rm (F_{III})}$, let $\alpha: Y\rightarrow \F_n$ be the composition of a sequence of blow-ups of $\F_n$ along $s$ points. Then $c_2(X)+4l$ is pseudo-effective.
\end{lem}

\begin{proof}
For the extremal ray $\R_+[l]$, we assume that $l$ is a smooth conic. It follows from (\ref{B:eqn11}) that $c_2(X)=(2s-4)l+f^*(-K_Y)(-K_X)+\Gamma$. Hence
\[ c_2(X)+4l=2sl+f^*(-(K_Y+\frac{1}{4}\Delta)).(-K_X)+\frac{1}{4}f^*\Delta.(-K_X)+\Gamma \]
is pseudo-effective, since both $-(4K_Y+\Delta)$ and $-K_X$ are nef.
\end{proof}

In general, when $X$ is a possibly singular quasi-projective variety, we can define the Chow ring $A(X)$ with its cap product, and for any coherent sheaf $\FF$ on $X$ with finite locally free resolution, we can define Chern classes $c_k(\FF)\in A^k(X)$ (cf.\ \cite{fu}). 

\begin{defn}\label{B:4.5}
Let $X$ be a terminal projective threefold, $S$ the singular locus of $X$ consisting of a finite number of points. Then $U=X\setminus S$ is a smooth quasi-projective threefold. We can give an alternative definition of $c_1(X),c_2(X)$ instead of using the general theory.
\begin{eqnarray}
c_1(X) & := & c_1(\TT_X|_U)\in A^1(X\setminus S)\tilde{\longrightarrow}A^1(X) \nonumber \\
c_2(X) & := & c_2(\TT_X|_U)\in A^2(X\setminus S)\tilde{\longrightarrow}A^2(X) \nonumber
\end{eqnarray}
\end{defn}

\begin{lem}\label{B:4.6}
With the notation as above. Let $\varphi: X'\rightarrow X$ be a resolution of $X$, such that $\varphi^{-1}$ is an isomorphism over $U$. Then $\varphi_*c_2(X')=c_2(X)$.
\end{lem}

\begin{proof}
Let $E$ be the exceptional locus of $\varphi$. Then there is an isomorphism $\varphi:V=X'\setminus E\rightarrow U=X\setminus S$. Thus $c_2(X)=c_2(\TT_X|_U)=\varphi_*c_2(\TT_{X'}|_V)$. Since $c_2(\TT_{X'}|_V)$ and $c_2(\TT_{X'})$ differ by a 1-cycle whose support is contained in $\Supp E$ and $\dim \varphi(E)=0$, we have $\varphi_*c_2(X')=c_2(X)$.
\end{proof}

\begin{lem}\label{B:4.7}
Assume that we are in case ${\rm (D_I)}$. Then $c_2(X)$ is pseudo-effective.
\end{lem}

\begin{proof}
In fact, $Y$ is a projective threefold with at most one terminal singular point. Let $r$ be the Gorenstein index of $Y$. Then $r=1$ or $2$ (cf.\ \cite{mo}). Since $-K_Y$ is nef and big, $c_2(Y)$ is pseudo-effective by Theorem \ref{B:2.2}(iii). It follows from Lemma \ref{B:4.6} that $f_*c_2(X)=c_2(Y)$. Let $H$ be any nef divisor on $X$. We may write that
\begin{eqnarray}
K_X & = & f^*K_Y+aE, \quad \hbox{where $a$ is positive rational}, \nonumber \\
H   & = & f^*f_*H+bE, \quad \hbox{where $b=H.E^2/E^3\leq 0$}. \nonumber
\end{eqnarray}
Then we have
\begin{eqnarray}
c_2(X).H & = & c_2(X).f^*f_*H+bc_2(X).E \nonumber \\
         & = & c_2(Y).f_*H+bc_2(X).E \nonumber \\
\chi(\OO_X) & = & \frac{1}{24}c_1(X).c_2(X) \nonumber \\
            & = & \frac{1}{24}(-f^*K_Y-aE).c_2(X) \nonumber \\
            & = & \frac{1}{24}c_1(Y).c_2(Y)-\frac{a}{24}c_2(X).E  \nonumber
\end{eqnarray}
It follows from the Singular Riemann-Roch formula (cf.\ \cite{re87}) that
\[ \chi(\OO_Y)=\frac{1}{24}c_1(Y).c_2(Y)+\frac{1}{24}(r-\frac{1}{r}) \]
Therefore $\chi(\OO_X)=\chi(\OO_Y)$ shows that $c_2(X).E=-(r-1/r)/a\leq 0$. It is easy to see that $f_*H$ is nef on $Y$, hence $c_2(Y).f_*H\geq 0$. Thus we have $c_2(X).H\geq 0$, which implies that $c_2(X)$ is pseudo-effective.
\end{proof}

\begin{lem}\label{B:4.8}
Assume that we are in case ${\rm (D_{II})}$. Furthermore, assume that $-K_Y$ is big. Then there exists a positive integer $n$ such that $c_2(X)+nl$ is pseudo-effective.
\end{lem}

\begin{proof}
We use the same notation as in Proposition \ref{B:3.9}. In general, we have
\[ c_2(X)=f^*c_2(Y)+C_1+(e-b)F\]
Note that $b-e$ is a non-negative integer and $c_2(Y)$ is pseudo-effective by assumption. If we exclude the case (C) in Proposition \ref{B:3.9}, then we obtain that $f_*H$ is nef on $Y$. This time, let $n=b-e, l=F$, then $c_2(X)+nl=f^*c_2(Y)+C_1$ is pseudo-effective by applying $H$ on each side.

We deal with the case (C) by the following lemma.
\end{proof}

\begin{lem}\label{B:4.9}
Assume that we are in one of the cases of Proposition \ref{B:3.9}. Furthermore, assume that $c_2(Y)$ is pseudo-effective. Then there exists a positive integer $n$ such that $c_2(X)+nl$ is pseudo-effective.
\end{lem}

\begin{proof}
Since $c_2(Y)$ is pseudo-effective, we may write $c_2(Y)=\lim_{k\to\infty}\xi_k$, where $\xi_k$ are effective 1-cycles. Let $a_k$ be the coefficient of $C$ in $\xi_k$ by writting $\xi_k=a_kC+R_k$. Then $\sup_k\{a_k\}<n$ for some suitable positive integer $n$. Given any nef divisor $H$ on $X$, $f_*H$ is nef on $Y$ possibly except along the curve $C$. By the same calculation as in Proposition \ref{B:3.9}, we have that $f_*H.C=H.f^*C=H.(C_1-rl)$, where $r=1$ for the cases (A) and (C), $r=2$ for the case (B). Let $s=0$ for the cases (A) and (B), $s=1$ for the case (C).

Thus we have
\begin{eqnarray}
c_2(X) & = & f^*c_2(Y)+C_1-sl \nonumber \\
(c_2(X)+2nl).H & = & f^*c_2(Y).H+(2n-s)l.H+C_1.H \nonumber \\
               & = & \lim_{k\to\infty}(a_kC+R_k).f_*H+(2n-s)l.H+C_1.H \nonumber \\
               & = & \lim_{k\to\infty}(a_kf^*C+(2n-s)l+f^*R_k).H+C_1.H\geq 0 \nonumber
\end{eqnarray}
which completes the proof.
\end{proof}

For case ${\rm (D_{III})}$, we only give an assumption denoted by ${\rm (AD_{III})}$.

\begin{asp}[${\bf AD_{III}}$]\label{B:4.10}
Assume that we are in case ${\rm (D_{III})}$. Then there exists a positive integer $n$ such that $c_2(X)+nl$ is pseudo-effective.
\end{asp}

\begin{prop}\label{B:4.11}
Let $X$ be a smooth projective threefold, $\R_+[l]$ an extremal ray on X. Assume that $c_2(X)+nl$ is pseudo-effective for some $n\in\N$. Then we can take an ample divisor $L$ on $X$, a sufficiently small $\varepsilon>0$, and the cone decomposition $\overline{NE}(X)=\R_+[l]+\sum_i\R_+[l_i]+\overline{NE}_\varepsilon(X)$ such that for any decomposition $c_2(X)=al+\sum_ib_il_i+z$, where $b_i\geq 0$, $z\in\overline{NE}_\varepsilon(X)$, we have $z.(-K_X)<1$.
\end{prop}

\begin{proof}
Let $H$ be a nef divisor on $X$, such that $H^\bot\cap\overline{NE}(X)=\R_+[l]$. It follows from \cite{mo} that $L=mH-K_X$ is ample for some integer $m\gg 0$. Fix $m$ and such ample divisor $L$, and take a $\varepsilon>0$ such that $\R_+[l]$ appears in the formula of the cone decomposition. Since $c_2(X)+nl$ is pseudo-effective, we have the following decomposition:
\[ c_2(X)=al+\sum_ib_il_i+z \]
where $b_i\geq 0$, $z\in\overline{NE}_\varepsilon(X)$. Since $H$ is nef and $l.H=0$, we have $z.H\leq c_2(X).H$. It follows from the definition of $\overline{NE}_\varepsilon(X)$ that
\begin{eqnarray}
z.(-K_X)\leq \varepsilon z.L=m\varepsilon(z.H)+\varepsilon z.(-K_X) \nonumber
\end{eqnarray}
Hence we have
\begin{eqnarray}
z.(-K_X)\leq \frac{m\varepsilon}{1-\varepsilon}z.H\leq \frac{m\varepsilon}{1-\varepsilon}c_2(X).H \nonumber
\end{eqnarray}
If $z.(-K_X)\leq 0$, then there is nothing to prove. Otherwise, $c_2(X).H>0$ and we may take $\varepsilon$ to be sufficiently small in advance so that $z.(-K_X)<1$.
\end{proof}

\begin{thm}\label{B:4.12}
Let $X$ be a smooth projective threefold such that $-K_X$ is nef, $\nu(-K_X)=2$ and $q(X)=0$. Assume that ${\rm (AD_{III})}$ holds. Then $c_2(X)$ is pseudo-effective.
\end{thm}

\begin{proof}
We use induction on the Picard number $\rho(X)$.

It is easy to see that $\rho(X)>1$. If $\rho(X)=2$, then only cases ${\rm (F_I),(F_{II}), (D_I)}$ and ${\rm (D_{II})}$ can occur. In cases ${\rm (F_I),(F_{II})}$ and ${\rm (D_I)}$, $c_2(X)$ is pseudo-effective. In case ${\rm (D_{II})}$, $-K_Y$ is ample since $\rho(Y)=1$. It follows from Lemma \ref{B:4.8} that $c_2(X)+nl$ is pseudo-effective for some positive integer $n$.

As in Proposition \ref{B:4.11}, we may take an ample divisor $L$ on $X$, a sufficiently small $\varepsilon>0$, and the cone decomposition
\[ \overline{NE}(X)=\R_+[l]+\sum_{i\in I}\R_+[l_i]+\overline{NE}_\varepsilon(X) \]
where the set of extremal rays $\{\R_+[l_i]\}_{i\in I}$ is fixed. If for some $\R_+[l_i]$, the corresponding extremal contraction is of type ${\rm (F_I),(F_{II})}$, or ${\rm (D_I)}$, then $c_2(X)$ is pseudo-effective. Otherwise all $\R_+[l_i]$ are of type ${\rm (D_{II})}$. Thus for each $i\in I$, there exists a positive integer $n_i$ such that $c_2(X)+n_il_i$ is pseudo-effective. Let $N$ be an integer such that $\max_{i\in I}\{n_i\}<N$.

Consider the decomposition of $c_2(X)$
\begin{eqnarray}
c_2(X)=al+\sum_{i\in I}b_il_i+z \label{B:eqn12}
\end{eqnarray}
where $a\geq -n$, $b_i\geq 0$ and $z\in\overline{NE}_\varepsilon(X)$. If $a\geq 0$, then $c_2(X)$ is pseudo-effective. So we may assume that $0>a\geq -n$.

Applying $(-K_X)$ to each side of (\ref{B:eqn12}), we have
\begin{eqnarray}
24 & = & c_1(X).c_2(X)=al.(-K_X)+\sum_{i\in I}b_il_i.(-K_X)+z.(-K_X) \nonumber \\
   & \leq & 4\sum_{i\in I}b_i+1 \nonumber
\end{eqnarray}
hence $\sum_{i\in I}b_i>5$, since $l_i.(-K_X)\leq 4$ and $z.(-K_X)<1$.

Consider the following pseudo-effective 1-cycle
\begin{eqnarray}
 &      & \sum_{i\in I}b_i(c_2(X)+n_il_i) \nonumber \\
 & \leq & (\sum_{i\in I}b_i)c_2(X)+N(\sum_{i\in I}b_il_i+z) \nonumber \\
 &  =   & (\sum_{i\in I}b_i+N)c_2(X)-Nal \nonumber
\end{eqnarray}
Since $\sum_{i\in I}b_i>5$, we have that $c_2(X)+\theta(-a)l$ is pseudo-effective, where $\theta=N/(N+5)<1$ is fixed. We may repeat the above argument to deduce that $c_2(X)+\theta^knl$ is pseudo-effective for any $k\in\N$. Thus $c_2(X)=\lim_{k\to\infty}(c_2(X)+\theta^knl)$ is pseudo-effective.

It follows from the above proof that if for every extremal ray $\R_+[l_i]$, there exists some positive integer $n_i$ such that $c_2(X)+n_il_i$ is pseudo-effective, then so is $c_2(X)$. This argument is used for the case $\rho(X)=2$ as well as in the inductive step.

Assume that when $\rho(X)=\rho-1\geq 2$, the conclusion holds. Let $\rho(X)=\rho$. Then cases ${\rm (F_I),(F_{II})}$ cannot occur. In case ${\rm (F_{III})}$, $c_2(X)+4l$ is pseudo-effective. In case ${\rm (D_I)}$, $c_2(X)$ is pseudo-effective. In case ${\rm (D_{II})}$, if $-K_Y$ is big, then $c_2(X)+nl$ is pseudo-effective. Otherwise if $-K_Y$ is nef but not big, then $c_2(Y)$ is pseudo-effective by the induction hypothesis or Theorem \ref{B:2.2}. It follows from Propositions \ref{B:3.9} and \ref{B:3.10} that $c_2(X)$ is pseudo-effective. In case ${\rm (D_{III})}$, we have that $c_2(X)+nl$ is pseudo-effective by ${\rm (AD_{III})}$. By the above argument, we can prove that $c_2(X)$ is pseudo-effective.
\end{proof}

\begin{rem}\label{B:4.13}
In case ${\rm (D_{III})}$, $-K_Y$ is no longer nef, just almost nef which is defined in \cite{ps}. So we cannot proceed by induction without ${\rm (AD_{III})}$. In order to get rid of ${\rm (AD_{III})}$, it is inevitable to extend the class of smooth projective threefolds with nef anticanonical divisors to the class of $\Q$-factorial terminal projective threefolds with almost nef anticanonical divisors to run the Minimal Model Program. I will deal with the terminal case in a subsequent paper. As a simple case, we show that ${\rm (AD_{III})}$ holds when $\rho(X)=3$ in the following proposition.
\end{rem}

\begin{prop}\label{B:4.14}
Assume that we are in case ${\rm (D_{III})}$. Furthermore, assume that $\rho(X)=3$. Then ${\rm (AD_{III})}$ holds.
\end{prop}

\begin{proof}
Let $f: X\rightarrow Y$ be an extremal contraction of type ${\rm (D_{III})}$ with $\rho(X)=3$. It is sufficient to prove that $c_2(Y)$ is pseudo-effective by Lemma \ref{B:4.9}. Since $\kappa(Y)=\kappa(X)=-\infty$, we can take an extremal contraction from $Y$. Note that $\rho(Y)=2$, so there are only four cases for such an extremal contraction $g:Y\rightarrow Z$ corresponding to an extremal ray $\R_+[l_0]$.

(\ref{B:4.14}.1) $Z\cong\PP^1$, and $g:Y\rightarrow Z$ is a del Pezzo fibration.

(\ref{B:4.14}.2) $Z\cong\PP^2$, and $g:Y\rightarrow Z$ is a conic bundle.

(\ref{B:4.14}.3) $Z$ is a terminal projective threefold with $-K_Z$ ample, and $g:Y\rightarrow Z$ is a birational morphism which contracts a divisor to a point.

(\ref{B:4.14}.4) $Z$ is a smooth projective threefold with $-K_Z$ ample, and $g:Y\rightarrow Z$ is a birational morphism which contracts a divisor to a smooth curve.

In the case (\ref{B:4.14}.1), the proof is similar to that of Lemma \ref{B:4.2}. Note that $M.(-K_Y)^2=(-K_Y|_M)^2\geq 0$ since $-K_Y$ is nef except along $C$ and $M$ is ample, and $(-K_Y)^3\leq 0$ by some direct computations. Hence we also have $b\geq 0$.

In the case (\ref{B:4.14}.2), the proof is similar to that of Lemma \ref{B:4.3}. Note that $F.(-K_Y)^2=(-K_Y|_F)^2\geq 0$.

In the case (\ref{B:4.14}.3), the proof is the same as that of Lemma \ref{B:4.7}.

In the case (\ref{B:4.14}.4), we only have that $c_2(Y)+n_0l_0$ is pseudo-effective for some positive integer $n_0$ by Lemma \ref{B:4.8}. Then by the same argument as in Theorem \ref{B:4.12}, we can prove that $c_2(Y)$ is pseudo-effective.
\end{proof}

\begin{cor}\label{B:4.15}
Let $X$ be a smooth projective threefold with $-K_X$ nef. If $\rho(X)\leq 3$, then $c_2(X)$ is pseudo-effective.
\end{cor}

\begin{proof}
This follows from Theorems \ref{B:2.2}, \ref{B:3.12} and \ref{B:4.12} and Proposition \ref{B:4.14}.
\end{proof}

\textsc{Graduate School of Mathematical Sciences, University of Tokyo, 3-8-1 Komaba, Meguro, Tokyo 153-8914, Japan}

\textit{E-mail address}: \texttt{xqh@ms.u-tokyo.ac.jp}

\end{document}